\documentclass[draft]{publmathdeb}
\usepackage{amsmath,amsfonts,amssymb}
\theoremstyle{plain}
\theoremstyle{definition}
\theoremstyle{remark}
\theoremstyle{plain}
\newtheorem{prop}{Proposition}
\newtheorem{thm}{Theorem}[section]

\newtheorem{coro}[thm]{Corollary}

\theoremstyle{definition}

\theoremstyle{remark}
\newtheorem{remark}{Remark}[section]

\theoremstyle{remark}

\theoremstyle{plain}

\begin{document}
\title{Irrational self-similar sets}
\author[K. Jiang]{Qi Jia, Yuanyuan Li and Kan Jiang}
\address[Q.Jia]{Department of Mathematics, Ningbo University,
People's Republic of China}
\email{2011071018@nbu.edu.cn}
\address[Y.Y.Li]{Department of Mathematics, Ningbo University,
People's Republic of China}
\email{1911071014@nbu.edu.cn}
\address[K. Jiang]{Department of Mathematics, Ningbo University,
People's Republic of China}
\email{jiangkan@nbu.edu.cn;kanjiangbunnik@yahoo.com}
\date{\today}
\begin{abstract}
Let $K\subset\mathbb{R}$ be a self-similar set defined on $\mathbb{R}$. It is easy to prove that  if the Lebesgue measure of $K$ is zero, then  for Lebesgue almost every $t$, $$K+t=\{x+t:x\in K\}$$ only  consists of  irrational or  transcendental numbers.  In this note, we shall consider some classes of self-similar sets, and explicitly construct such $t$'s. Our main idea is from the $q$-expansions.
\end{abstract}
\submitted{February 2, 2021}
\maketitle
\section{Introduction}
Given $m\in \mathbb{N}_{\geq 2}$,
let  $\{\phi_{i}\}_{i=1}^{m}$ be an iterated function system of contractive similitudes on $\mathbb{R}$ defined as
 \begin{equation}\label{Hutchinson formula}
 \phi_{i}(x)=\lambda_{i}x+b_{i},\;\; i=1,\ldots ,m,
\end{equation}
where $0<|\lambda_{i}|<1$ and  $b_{i}\in \mathbb{R}.$ Hutchinson  \cite{Hutchinson} proved that there exists a unique non-empty compact set, denoted by $K$, such that
$$K=\cup_{i=1}^{m}\phi_i(K).$$
We call $K$ the self-similar set or attractor with respect to the IFS
$\{\phi_{j}\}_{j=1}^{m}$, see \cite{Hutchinson} for further details.
In particular, some  $q$-adic expansions can be generated by some IFS.  Let $q\in \mathbb{N}_{\geq 3}$. Define  a digit set $$ \mathcal{A}\subset \{0,1,2,\cdots,q-1\}.$$
Without loss of generality, we  assume $\mathcal{A}$  contains at least two integers.
Then the following set
\begin{equation}\label{IFS}
K=\left\{x=\sum_{i=1}^{\infty}\dfrac{a_i}{q^i}:a_i\in \mathcal{A} \right\}.
\end{equation}
is clearly a self-similar set.
Its IFS is $$\left\{\phi_i(x)=\dfrac{x+c_i}{q},c_i\in \mathcal{A}\right\} .$$
In other words, for any $x\in K$, we can find a  sequence $(a_i)\in \mathcal{A}^{\mathbb{N}}$ such that
$$x=\sum_{i=1}^{\infty}\dfrac{a_i}{q^i}.$$ We call $(a_i)$ a $q$-adic expansion or coding of $x$.  Evidently, if $\mathcal{A}$ does not contain consecutive digits, then any point in $K$ has a unique $q$-adic expansion.

Let $K\subset\mathbb{R}$ be a self-similar set. Usually, $K$ simultaneously  contains rational and irrational numbers.
It is natural to ask whether $K$ only consists of irrational numbers. For simplicity,  we call a set $E\subset \mathbb{R}$ an irrational set if all the numbers in $E$ are irrational. We say that $E$ is a transcendental set if all the numbers in $E$ are transcendental \cite{Allouche}. We first give a  result on the  existence of  $t$ such that   $K+t$ is an irrational set or a  transcendental set.  We,  in fact, can prove more general results. Namely, the following set
$$\{t\in \mathbb{R}: K+t \mbox{ is an  irrational set or  a transcendental set} \}$$ is large   from the measure theoretical or topological perspective.
\begin{prop}\label{key1}
If  $E\subset \mathbb{R}$ is a set with zero Lebesgue measure and $D\subset\mathbb{R}$ is a countable set. Then for Lebesgue almost every $t$, $$(E+t)\cap D =\emptyset.$$
Moreover, if $E\subset \mathbb{R}$ is a set of  first category, then  except a set of the first category with respect to the parameter $t$,  we have that
$$(E+t)\cap D=\emptyset.$$
\end{prop}
Proposition \ref{key1}   implies many results. We only list the following two consequences. The reader may find other  similar statements. For instance, analogous results can be obtained from the topological perspective.
\begin{coro}\label{cor1}
Let $\{E_i\}_{i=1}^{\infty}\subset \mathbb{R}$ be a set sequence such that each $E_i, i\in \mathbb{N}^{+}$  is a  Lebesgue  null set. Then for  Lebesgue almost every $t$, we have that
$$\cup_{i=1}^{\infty}(E_i+t) $$ is  an  irrational set or  a  transcendental set.
\end{coro}
\begin{coro}\label{cor2}
 Let $E$ be a set with zero Lebesgue measure. Then for   Lebesgue almost every $t$,
$$E+it=\{e+it:e\in E\},\dfrac{E}{it}=\left\{\dfrac{e}{it}:e\in E\setminus\{0\}\right\}, i\in \mathbb{Z}\setminus \{0\}$$ are  simultaneously irrational sets  or  transcendental sets.
Moreover, for  Lebesgue almost every $t$,
$$E+t, E-t, tE, \dfrac{E}{t}$$are  simultaneously irrational sets  or  transcendental sets.
Here $$tE=\{te, e\in E\setminus\{0\}\}.$$
\end{coro}
Although Proposition \ref{key1} gives the existence of $t$, usually it is not easy to find a such $t$ explicitly.  The main aim of this paper is to consider some self-similar set $K$, and construct explicitly $t$
such that $K+t$ is an irrational set.
Before we introduce the first concrete example, we define a Liouville number.
Let $q\in \mathbb{N}_{\geq 3}$. We define
$$s(q):=\sum_{n=1}^{\infty}\dfrac{1}{q^{n!}}.$$
In what follows, we shall use this number $s(q)$.  Sometimes, we may simply denote it by $s.$
\begin{thm}\label{Main}
Let $K$ be the self-similar set generated by
$$\left\{f_1(x)=\dfrac{x}{q},f_2(x)=\dfrac{x+q-1}{q}\right\}$$ with  $  q\in \mathbb{N}_{\geq 3}.$
Then for any $\ell\in \{1,2,\cdots, q-2\}$ we have
$$K-\ell s=\{x-\ell s:x\in K\}\subset \mathbb{Q}^c,$$
In particular, when $q=3$, $K$ is the middle-third Cantor set (denoted by $C$).  Therefore, we have
$$C- s(3)  \subset \mathbb{Q}^c.$$
\end{thm}
\begin{remark}
\begin{itemize}
\item[(1)] Our result can be strengthened, i.e. for any $r\in \mathbb{Q}$, we have
$$K-\ell s-r\subset \mathbb{Q}^c.$$
It is easy to see that $K-(q-1)s$ may contain some rational numbers.  For instance, we let $q=3$. We claim that $C-2s$ may contain $0$. This is because $$x=\sum_{n=1}^{\infty}\dfrac{2}{3^{n!}}\in C$$ and $x-2s=0.$
 We do not know whether each  $K- \ell s, 1\leq \ell\leq q-2$ is a transcendental set. However,  we can conclude that for any $1\leq \ell\leq  q-2$, all the numbers in  $K- \ell s$ cannot only contain the Liouville numbers as the set of all the Liouville numbers  has Hausdorff dimension zero whilst the translation of $K$ has Hausdorff dimension $$\dfrac{\log 2}{\log q}.$$
\item[(2)] The Liouville number $s$ is not the unique way which guarantees $$K-\ell s\subset \mathbb{Q}^c.$$ The number $s$ can be constructed generally  as follows.

    Let $g:\mathbb{N}^{+}\to \mathbb{N}^{+}$ be an integer function such that $$\lim_{n\to \infty}(g(n+1)-g(n))=+\infty.$$
    Then we define
    $$s^{*}=\sum_{n=1}^{\infty}\dfrac{1}{q^{g(n)}}.$$
    In particular, we may let $g(n)=n!, n^n, 2021n^3+2020n^2+2019n, \cdots. $
    Under these constructions, we still have the desired result in the above theorem.  In what follows, we always let $g(n)=n!$. Some other constructions are allowed.

    \item[(3)] Similar results can be obtained for the $p$-adic numbers ($p$ is some prime number). As in the setting of $p$-adic system, we still have that
    a number in $\mathbb{Q}_p$ is rational if and only if it has an eventually periodic $p$-adic expansion. The proof is almost the same as our idea. We leave this to the reader.
    \item[(4)] If we want to obtain the following result, i.e. for any $i\in \mathbb{Z}\setminus\{0\}$,
    $$K+is^{\prime}\subset \mathbb{Q}^c $$
 then we may let $s^{\prime}$ be a normal number, see \cite{Jiang} for more general result.  For the Liouville number, it is not enough. A counterexample is constructed in the above Remark (1).

\end{itemize}

\end{remark}
With similar discussion, we can prove the following result.
\begin{thm}\label{Main1}
Let $K_1$ be the attractor of the following IFS
$$\left\{f_i(x)=\dfrac{x+a_i}{q}\right\},$$ where
$$ a_i\in \mathcal{A}\subset \{0,1,2,\cdots,q-1\},  q\in \mathbb{N}_{\geq 3}.$$
Suppose
\begin{itemize}
\item [(1)] $1\notin \mathcal{A}, q-1\in \mathcal{A}$;
\item [(2)]   $\mathcal{A}$ contains no consecutive integers.
\end{itemize}
Then  we have
$$K_1-s=\{x-s:x\in K_1\}\subset \mathbb{Q}^c,$$
\end{thm}
\begin{remark}\label{remark}
The proofs of Theorems \ref{Main} and \ref{Main1} are motivated by some classical techniques in the setting of $q$-expansions. In particular, we are motivated by  some carries occurring  in the Fibonacci base \cite{GS}, i.e. a word $100$ can be replaced by the word $011$ in base $\dfrac{\sqrt{5}+1}{2}$ .  The reader may find many useful ideas in the references \cite{RSK,KarmaMartijn,MK,EJK,GS,SN}.
\end{remark}
Theorem \ref{Main1} has the following corollary.
\begin{coro}
Let $J$ be the attractor of the IFS
$$f_i(x,y)=\left(\dfrac{x+c_i}{q}, \dfrac{x+d_i}{q}\right),$$ where
$(c_i,d_i)\in \mathcal{A}_1\times \mathcal{A}_2$, and $\mathcal{A}_i, i=1,2,$ satisfies the conditions in Theorem \ref{Main1}.
Then $$J+\left(-s,-s\right)=\left\{\left(x-s, y-s\right):(x,y)\in J\right\}\subset \mathbb{Q}^c \times \mathbb{Q}^c.$$
\end{coro}
We may  construct a self-similar set $J_1\subset \mathbb{R}^2$ with zero Lebesgue measure such that for any $(x,y)\in \mathbb{R}^2$, we always have
$$J_1+(x,y)\nsubseteqq \mathbb{Q}^c  \times \mathbb{Q}^c .$$ The key point is that a two dimensional set with zero Lebesgue measure may still contain some line segment.

For different self-similar sets, we may give a uniform translation such that
the translation sets are irrational sets.
\begin{thm}\label{Main2}
Let $J_2$ and $J_3$ be the attractors of the IFS's
$$\left\{f_1(x)=\dfrac{x}{q^n}, f_2(x)=\dfrac{x+q^n-1}{q^n}\right\}$$
and
$$\left\{g_1(x)=\dfrac{x}{q^m}, g_2(x)=\dfrac{x+q^m-1}{q^m}\right\},$$ respectively, where $q\in \mathbb{N}_{\geq 3}, n,m\in \mathbb{N}_{\geq 2}.$ Then
$$J_2+s\subset \mathbb{Q}^c, J_3+s\subset \mathbb{Q}^c,$$
where
$$s=-\sum_{k=1}^{\infty}\dfrac{1}{q^{(mn)k!}}.$$
\end{thm}
\begin{remark}
We can give more similar results. For instance, in Theorem \ref{Main2}, we are allowed to consider a general digit set as defined in Theorem \ref{Main1}. Moreover, in Theorem  \ref{Main1}, we may investigate consecutive translations as Theorem \ref{Main}. In other words, we are able to find some results for a general  digit set, and to analyze the consecutive translations of some fractal sets. We leave these generalizations to the reader.
\end{remark}
This paper is organized as follows. In Section 2, we give  the  proofs of the main results. In Secion 3, we pose some problems.
\section{Proofs of the main results}
\begin{proof}[\textbf{Proof of Proposition \ref{key1}}]
We only prove   the first statement as the second proof is similar.
Suppose that for some $t$,  $E+t$ contains some number in $D$, i.e.
$$D\cap (E+t)\neq \emptyset.$$ Therefore,  there exists some $a\in D, e\in E $ such that
$$-t=e-a\in E-a.$$ Subsequently,
$$-t\in \cup_{a\in D}(E-a),$$ which yields that
$$m(\cup_{a\in D}(E-a))\leq \sum_{a\in D}m(E-a)=0, $$
where $m(\cdot)$ denotes the Lebesgue measure.
\end{proof}
The proofs of Corollaries \ref{cor1}  and \ref{cor2} are similar to Proposition \ref{key1}. We leave them to the reader.
\begin{proof}[\textbf{Proof of Theorem \ref{Main}}]
Note that for any $x\in K$ there exists a unique sequence $(x_i)\in\{0,q-1\}^{\mathbb{N}}$ such that
$$x=\sum_{i=1}^{\infty}\dfrac{x_i}{q^i}=:(x_1x_2\cdots )_q,$$ and
the infinite sequence $(x_i) $ is called a (unique) $q$-adic coding of $x$. Define
    $$s=\sum_{n=1}^{\infty}\dfrac{1}{q^{n!}}=(s_1s_2\cdots)_q,$$
    where
  \begin{equation}\label{1}
  s_i=\left\{
\begin{array}{ll}
1 & \text{if } i=k! \text{ for some }k\in \mathbb{Z}^{+} \\
0 & \text{otherwise} .
\end{array}%
\right.
  \end{equation}
 Take $x\in K$ with its unique coding $(x_i)\in\{0,q-1\}^{\mathbb{N}}, $ and take $\ell \in \{1,\cdots, q-2\}$. It suffices to prove that
  \begin{equation}\label{2}
  x-\ell s=\sum_{i=1}^{\infty}\dfrac{x_i-\ell s_i}{q^i}\notin \mathbb{Q}.
  \end{equation}
  Write
  $a_i=x_i-\ell s_i$ for $i\geq 1$. Then $$a_i\in \{0,q-1\}-\{0,\ell\}=\{-\ell,0,q-1-\ell,q-1\}.$$ Furthermore, by (\ref{1}) it follows that $a_i\in \{-\ell, q-1-\ell\}$ if $i=k!$ for some  $k\in \mathbb{N},$ and otherwise $a_i\in \{0,q-1\}.$ It is well-known that a number $y\in \mathbb{R}$ is rational if and only if it has an eventually periodic $q$-adic coding in $\{0,1,\cdots, q-1\}^{\mathbb{N}}$. Observe that the sequence $(a_i)$ may have negative digit $-\ell$. Therefore, our strategy to prove (\ref{2}) is that we first construct a $q$-adic coding $(b_i)\in \{0,1,\cdots, q-1\}^{\mathbb{N}}$ of $x-\ell s$ based on the sequence $(a_i)$, and then show that the coding $(b_i)$ is not eventually periodic.
    Without loss of generality we may assume that $$x=(x_1x_2\cdots )_q\in K$$ with $x_1=q-1$.
  This is due to the similarity structure of $K$, i.e. for any $y\in [0,1/q]\cap K$, then there exists some $x\in [1-1/q,1]\cap K$ such that
  $$x=y+1-1/q.$$
  Therefore, $$x\notin \mathbb{Q}\Leftrightarrow y\notin \mathbb{Q}.$$

  In the following we  may assume that $(x_i)$ contains infinitely many zeros and infinitely many $(q-1)$'s. If $(x_i)$ ends with $(q-1)^{\infty}$ or   $0^{\infty}$ (in what follows, we denote $a^k$ by $k$ consecutive concatenation of the digit $a$), then clearly $x$ is rational. So $x-\ell s\notin \mathbb{Q}$ as $s$ is a transcendental number. In other words, we finish the proof.

  Note that $a_1\in \{q-1-\ell, q-1\}.$ If $a_{2!}=-\ell$, then we replace $a_1a_2$ by $b_1a^{\prime}_1=(a_1-1)(q-\ell).$ Now we look at the next place where $a_{(k+1)!}=-\ell$ for some smallest $k\geq 2$.    If there exists a largest index $i\in \{k!+1, \cdots, (k+1)!-1\}$ such that $a_i=q-1$, then we replace the word
  $$a_{k!}^{\prime}a_{k!+1}\cdots a_{(k+1)!}$$ by
  $$b_{k!}\cdots b_{(k+1)!-1} a_{(k+1)!}^{\prime}=a_{k!}^{\prime} a_{k!+1}\cdots a_{i-1}(a_i-1)(q-1)^{(k+1)!-i-1}(q-\ell).$$
  Otherwise,  $a_i=0$ for all $k!<i<(k+1)!$. Note that by our previous replacement, we have  $a_{k!}^{\prime}>0.$ Then we replace the word
    $$a_{k!}^{\prime}a_{k!+1}\cdots a_{(k+1)!}$$ by
  $$b_{k!}\cdots b_{(k+1)!-1} a_{(k+1)!}^{\prime}=(a_{k!}^{\prime}-1) (q-1)^{(k+1)!-k!-1}(q-\ell).$$
  Proceeding this argument indefinitely we obtain a $q$-adic coding $(b_i)\in \{0,1,\cdots, q-1\}^{\mathbb{N}}$ of $x-\ell s.$

  In the following we will prove by contradiction that $(b_i)$ is not eventually periodic. Suppose on the contrary that $$(b_i)=b_1\cdots b_N(b_{N+1}\cdots b_{N+p})^{\infty}$$ for some $N,p\in \mathbb{N}.$ Then
   \begin{equation}\label{3}
   b_{i+p}=b_i,\forall i\geq N.
   \end{equation}
   Take $k\geq \max\{p, N\}$. Observe that $b_{k!}\in\{q-\ell-2, q-\ell-1,q-\ell\}$. We split the proof into the following three cases.

    Case \uppercase\expandafter{\romannumeral1}. $b_{k!}=q-\ell-2$, Then $b_{k!+1}\cdots b_{(k+1)!-1}=(q-1)^{(k+1)!-k!-1}$.This, together with (\ref{3}), implies that $(b_i)$ ends with $(q-1)^{\infty}$. However, by our construction the sequence $(b_i)$ cannot end with $(q-1)^{\infty}$, leading to a contradiction.

      Case \uppercase\expandafter{\romannumeral2}. $b_{k!}=q-\ell-1$. Then $1\leq q-\ell-1\leq q-2.$ By (\ref{3}) it follows that there exist more than one index $i\in \{k!+1, \cdots, (k+1)!-1\}$ such that $b_i=q-\ell-1.$ This is again leads to a contradiction with our construction of $(b_i). $

       Case \uppercase\expandafter{\romannumeral3}.   $b_{k!}=q-\ell$.  Then there exists $i\in \{k!+1, \cdots, (k+1)!-1\}$ such that
       $$b_i\cdots b_{(k+1)!-1}=(q-2)(q-1)^{(k+1)!-i-1}.$$ If $(k+1)!-i-1\geq p$, then we can conclude by (\ref{3}) that $(b_i)$ ends with $(q-1)^{\infty}$. Otherwise, we can conclude that there exist more than one index
       $i\in \{(k+1)!+1, \cdots, (k+2)!-1\}$ such that $b_i=q-2.$ However, both cases will lead to a contradiction with our construction of $(b_i).$

       Therefore, by Cases \uppercase\expandafter{\romannumeral1}-\uppercase\expandafter{\romannumeral3} it follows that $(b_i)$ is not eventually periodic, and thus $x-\ell s=(b_1b_2\cdots )_q\notin \mathbb{Q}. $ This completes the proof.
\end{proof}
\begin{proof}[\textbf{Proof of Theorem \ref{Main1}}]
The proof of Theorem \ref{Main1} is similar to that of  Theorem \ref{Main}.
It  suffices to  prove that  for any
$$x\in K\cap [1-1/q,1] $$ we have
$$r=x-s\notin  \mathbb{Q}.$$
Let $$x=\sum_{i=1}^{\infty}\dfrac{x_i}{q^i}$$ with $(x_i)\in \mathcal{A}^{\mathbb{N}}, x_1=q-1$. Let $(s_i)\in \{0,1\}^{\mathbb{N}}$ with $s_i=1$ iff $i=k!.$
Then it is easy to check that
$$r=\sum_{i=1}^{\infty}\dfrac{x_i- s_i}{q^i}\in (0,1),$$ and that
$$x_1- s_1=q-1-1\geq 1, x_i-s_i\in\{-1,0,1,2,\cdots, q-1\}, i\geq 2.$$
We denote $b_i:=x_i-s_i. $ Therefore, we have
$$r=\sum_{i=1}^{\infty}\dfrac{b_i}{q^i}.$$
Then we may implement the idea from Theorem \ref{Main} and
rearrange the expansion $(b_i)$. Without loss of generality, we still use $(b_i)$ to denote the new  rearranged expansion. We can prove
the following two statements.
 \begin{itemize}
 \item[(1)]
  If $i=k!$ for some  large $k\geq 2$, then
    $$b_i\in  \{(\mathcal{A}-1)\cap \{0,1,\cdots, q-1\}\}\subset \mathcal{A}^c.$$
 \item[(2)]
If there exists some  large $k\geq 3$ such that $k!<i<(k+1)!$, then $$b_i\in \mathcal{A}.$$
  \end{itemize}
  Therefore, the expansion $(b_i)$ is not eventually periodic, which yields that $r\notin \mathbb{Q}.$
\end{proof}
\begin{proof}[\textbf{Proof of Theorem \ref{Main2}}]
The proof is similar to that of Theorem \ref{Main}. We only give an outline.
It suffices to consider  the numbers in
$$ \left[\dfrac{q^n-1}{q^n},1\right]\cap J_2,  \left[\dfrac{q^m-1}{q^m},1\right]\cap J_3.$$
More specifically, we let $r_1=x_1+s, r_2=x_2+s$, where
$$x_1\in \left[\dfrac{q^n-1}{q^n},1\right]\cap J_2, x_2\in \left[\dfrac{q^m-1}{q^m},1\right]\cap J_3.$$
Note that $r_1, r_2\in (0,1)$. Therefore, both of them have some $q$-expansions.  We denote their expansions by
$$(u_i)\in \{-1,0,q^n-1, q^n-2\}^{\mathbb{N}}, (v_i)\in \{-1,0,q^m-1, q^m-2\}^{\mathbb{N}},$$
respectively. By means of the idea we used in Theorem \ref{Main}, we rearrange these two expansions. Again, we still use $(u_i)$ and $(v_i)$ to denote these expansions.
We can prove the following statements:
\begin{itemize}
\item[(1)] if $i$ is sufficient large, then $u_i=q^n-2$, where $i=(mn)k!$ for some large $k\geq 1;$
\item[(2)] if $i$ is sufficient large,  and  $(mn)k!<i<(mn)(k+1)!$ for some large $k$, then $u_i\in \{0,q^n-1\}$;
\item[(3)] if $i$ is sufficient large, then $v_i=q^m-2$, where $i=(mn)k!$ for some large $k\geq 1;$
\item[(4)] if $i$ is sufficient large,  and  $(mn)k!<i<(mn)(k+1)!$ for some large $k$, then $v_i\in \{0,q^m-1\}$.
\end{itemize}
Clearly, under these conditions, $(u_i)$ and $(v_i)$ are not eventually periodic. Therefore, we prove that
$r_1$ and $r_2$ are not rational.
\end{proof}
\section{Final remark}
There are many problems can be asked. We list the following questions.
\begin{itemize}
\item [(1)] Can we find a uniform $t$ such that $$C_{1/3}+t, C_{1/4}+t\subset \mathbb{Q}^c,$$ where $C_{1/3}$ and $C_{1/4}$ are the middle-third and  middle-$1/2$ Cantor sets, respectively. By Corollary \ref{cor1}, we can find in theory many such $t$. We conjecture that
    $$t=-\sum_{n=1}^{\infty}\dfrac{1}{3^{n!}}-\sum_{n=1}^{\infty}\dfrac{1}{4^{n!}}$$ may work.
\item [(2)] How can we find some explicit   self-similar set which contains only transcendental numbers.

    \item [(3)] We do not know whether $C_{1/3}+e$ or $C_{1/3}+\pi$ is a transcendental  set.

    \item [(4)] For any $1<q<2$,  how can we find an explicit $t$ such that $U_q+t$ is a transcendental set, where $U_q$ denotes the univoque set \cite{RSK,KarmaMartijn,MK,EJK,GS}.
\end{itemize}

\section*{Acknowledgements}
We are grateful to the anonymous referees for many suggestions and comments.
This work is  supported by the National Natural Science Foundation of China (No.11701302), and by the Zhejiang Provincial Natural Science Foundation of China with
No.LY20A010009. The work is
also supported by K.C. Wong Magna Fund in Ningbo University.

\end{document}